\documentclass[11pt]{amsart}

 \newtheorem{thm}{Theorem}
\newtheorem{lem}{Lemma}[section] 
\newtheorem{prop}[thm]{Proposition} 
\newtheorem*{thm*}{Theorem}

\newcommand{\R}{\mathbb{R}} \newcommand{\Q}{\mathbb{Q}}
\newcommand{\N}{\mathbb{N}} 
\newcommand{\Z}{\mathbb{Z}} \newcommand{\T}{\mathbb{T}}
\newcommand{\qn}{q_{s(n)}} \newcommand{\qnn}{q_{s(n)+1}}

\newcommand{\norm}[1]{\bigl\| #1 \bigr\|}
\renewcommand{\norm}[1]{|\!|\!| #1 |\!|\!|}

\newcommand{\var}{\operatorname{var}}

\begin{document}
\newlength{\hypwidth}
\settowidth{\hypwidth}{[H3]}

\begin{abstract}
  We provide sufficient conditions on a positive function so that its
  associated special flow over any irrational rotation is either weak
  mixing or $L^2$-conjugate to a suspension flow. 
\end{abstract}    
    
\title[A Spectral Dichotomy For Some Special Flows]{A Dichotomy
  between Discrete and Continuous Spectrum for a Class of Special
  Flows over Rotations.}

\author{B. Fayad}
\thanks{Prepared during B. Fayad's visit to the Pennsylvania
     State University, Fall 2001}
\address{LAGA, Universite Paris 13, CNRS UMR 7539.}
\email{fayadb@math.univ-paris13.fr}
\author{A. Windsor}
\thanks{A. Windsor partially supported by NSF grant DMS 0071339}
\address{Pennsylvania State University, University Park, PA 16802, USA}
\curraddr{University of Manchester, Oxford Road, Manchester, M13 9PL, UK}
\email{awindsor@maths.man.ac.uk}

 \thanks{The authors gratefully acknowledge Y. Katznelson for
   suggesting to us the use of Salem and Zygmund's paper
   \cite{MR9:181d}, and A. Katok for suggesting the problem to us and
   for helpful conversations. }


\maketitle

\section{Introduction}

In his I.C.M. address of 1954 \cite{MR16:36g}, Kolmogorov raised a
number of questions concerning reparameterization of irrational linear flows
on $\mathbb{T}^n$, or equivalently special flows over translations.
One of them was to determine what kind of spectral properties could be
displayed by the unitary operator associated to the special flow built
over an irrational rotation on the circle and under an analytic roof
function.

Kolmogorov noticed that if the rotation angle $\alpha$ is not very
well approximated by rational angles, e.g. $\alpha$ Diophantine, and
if the roof function $\varphi$ is a strictly positive real analytic
function, then the special flow $T^t_{(\alpha, \varphi)}$ built over
the rotation $R_\alpha$ and under the function $\varphi$ is
analytically conjugate to a constant time suspension over $R_\alpha$,
i.e.  to an irrational linear flow on $\T^2$. The argument, based on
solving an additive cohomological equation, also proves that, for any
irrational angle $\alpha$, if the roof function is a strictly positive
trigonometric polynomial then the special flow $T^t_{(\alpha,
  \varphi)}$ is analytically conjugate to an irrational linear flow on
$\T^2$.

Later, \v{S}hklover proved that for any strictly positive real analytic
function that is not a trigonometric polynomial, there exists an
irrational angle $\alpha$ such that the special flow
$T^t_{(\alpha,\varphi)}$ has continuous spectrum \cite{MR37:1737}.  Thus,
for analytic functions $\varphi$ that are not trigonometric
polynomials, both continuous and discrete spectra can be obtained
depending on $\alpha$.

These are not the only possibilities. In a recent work, the authors
together with A. Katok, have proved that for every Liouvillean angle
$\alpha$ there exists a strictly positive $C^\infty$ function
$\varphi$ such that the special flow $T^t_{(\alpha, \varphi)}$ has
mixed spectrum.  When the angle $\alpha$ is exceptionally well
approximated by rational numbers the function $\varphi$ can be made
analytic \cite{MR1901073}. The roof functions constructed in
\cite{MR1901073} have blocks of relatively large Fourier coefficients
which appear in a lacunary progression.  The possibility that mixed
spectra would be precluded for roof functions with more regular decay
of Fourier coefficients was raised in \cite{MR1858535}, a slightly
reworked part of earlier unpublished notes \cite{Unpublished}. In
these notes the authors observe that for a function such as
$$ 
\varphi(x) = \sum_{n \in \Z} 2^{-|n|} \cos ( 2 \pi n x)
$$ the special flow $T^t_{(\alpha, \varphi)}$ is conjugated to a
linear flow if $\alpha$ is such that there exists a $c > 0$ so that
for all $p \in \Z$ and $q \in \N$,
$$2^q \bigl| \alpha - \frac{p}{q}\bigr|>c.
$$ Conversely, they show that a sufficient condition for weak mixing
is the existence of sequences $\{p_n\}$ and $\{q_n\}$ such that
$$2^{q_n} q_n \bigl| \alpha - \frac{p_n}{q_n} \bigr| \rightarrow 0.$$
To prove weak mixing they use a criterion involving the distribution
of the Birkhoff sums of the roof function $\varphi$,
$$S_m \varphi(x) = \sum_{k=0}^{m-1} \varphi( x + k \alpha),$$
along a sequence $m_n$ satisfying $R_\alpha^{m_n} \rightarrow
Id$. They are able to choose a sequence where each $m_n$ is a multiple
of a single frequency $q_n$.  

In order to bridge the gap between the conditions above and prove a
full dichotomy depending on $\alpha$ between continuous and discrete
spectrum for $T^t_{(\alpha,\varphi)}$, we consider the distribution of
the Birkhoff sums $S_m \varphi(x)$ along sequences $m_n$ which again
have the property that $ R^{m_n}_\alpha \rightarrow Id$ but which
involve multiple frequencies $q_n$. We use conditions on the
regularity of the decay of the Fourier coefficients (c.f. [H1]) to
extract for each $\alpha$ a lacunary representative of the additive
cohomology class of $\varphi$.  If $\varphi$ is not an $L^2$
coboundary, we use the central limit theorem for lacunary series to
study the distribution of the Birkhoff sums of its lacunary
representative and prove weak mixing. Our motivation for examining
lacunary series was a result of M. Herman \cite{Herman}:
\begin{thm*}
  If $\varphi$ is a lacunary Fourier series, and the equation
  $$\psi(x+ \alpha) - \psi(x) = \varphi(x)$$
  has a measurable solution
  $\psi$, then in fact the equation has a solution in $L^2$.
\end{thm*}
Our result can be viewed as a rigidity result that covers the
multiplicative equation too. Either the additive equation admits an
$L^2$ solution or else the multiplicative equation admits no solution.

\section{A Weak Mixing Dichotomy for Special Flows}

\begin{thm}\label{thm:Dichotomy}
  Let $\varphi : \T \rightarrow \R^+$ be a $C^{3}$ function given by
  \begin{displaymath}
    \varphi (x) = \sum_{m \in \Z} c_m e^{2 \pi i m x},
  \end{displaymath}
  where the coefficients satisfy the regularity conditions
  \newcounter{hypcounter}
  \begin{list}{\rm[H\arabic{hypcounter}]}{\usecounter{hypcounter}
      \setlength{\labelsep}{2\labelsep}
      \setlength{\labelwidth}{\hypwidth}
      \addtolength{\labelwidth}{\labelsep}}
  \item there exist $C_m$ such that
    $$\sum_{l=2}^\infty |c_{l m} |^2 \leq C_m |c_m|^2
    $$
    and 
    $$\sum_{m=1}^\infty C_m < \infty.
    $$
    
  \item there exists $0< K_1 < 1/4$ such that 
    $$\sum_{l=2}^\infty |c_{l m}| < K_1 |c_m|
    $$
    for all $m$ sufficiently large,
  \item there exists $K_2 > 0$ such that 
    $$\sum_{l=2}^\infty |l c_{l m}| < K_2 |c_m|
    $$
    for all $m$ sufficiently large.

  \end{list}

  Then for all $\alpha \in \R \backslash \Q$ we have either
  \begin{enumerate}
  \item the special flow $T^t_{(\alpha, \varphi)}$ is weak mixing, or

  \item the special flow $T^t_{(\alpha, \varphi)}$ is $L^2$ conjugate
  to a suspension flow.
  \end{enumerate}
\end{thm}

\subsection{Examples:} The hypotheses [H1], [H2], and [H3] restrict
  the coefficients along arithmetic progressions. Thus relatively
  prime frequencies do not influence each other directly.

  \begin{lem}
    A positive $C^3$ function $\varphi$ given by
    \begin{displaymath}
      \varphi(x)  = c_0 + \sum_{ |p| prime} c_p e^{2 \pi i p x}
    \end{displaymath}
    satisfies {\rm [H1]}, {\rm [H2]}, and {\rm [H3]}.
  \end{lem}
  
  Regular exponential decay along the appropriate arithmetic
  progressions will also suffice.

  \begin{lem}
    A positive function $\varphi$ given by
    \begin{displaymath}
      \varphi (x) = \sum_{m \in \Z} c_m e^{2 \pi i m x},
    \end{displaymath}
    where the coefficients satisfy the regularity condition 
    \begin{displaymath}
      C_1 e^{- k_1 |m|} \leq |c_m| \leq C_2 e^{-k_2 |m|}
    \end{displaymath}
    with $1 \leq k_1 / k_2 < 2$, satisfies {\rm [H1]}, {\rm [H2]}, and
    {\rm [H3]}.
  \end{lem}

\section{The Tools}

\subsection{Arithmetic}\label{sec:Arithmetic}

Associated to each $\alpha \in \R \backslash \Q$ there is an infinite
sequence of natural numbers $\{ q_n \}$ which we call the sequence of
best returns. This sequence can be computed using continued fractions
as the denominators of successive convergents. We introduce the
notation
\begin{displaymath}
  \norm{x} = \inf_{p \in \Z} | x - p|
\end{displaymath}
to measure the distance of $x$ from $0$ in $\T$. We call $\norm{ q
  \alpha}$ the quality of the return $q$. We have
\begin{displaymath}
  \norm{q_n \alpha} < \norm{q \alpha} 
\end{displaymath}
for $1 \leq q < q_n$ and $q_n < q < q_{n+1}$. This justifies our
best return nomenclature. 

We will use two lemmas from the theory of continued fractions, see
\cite{MR98c:11008}. The first relates the speed of growth of the best
returns $\{ q_n \}$ with the quality of returns.
\begin{lem}\label{lem:ContinuedFractionsEstimates}
  Let $\alpha \in \R^+ \backslash \Q$ and let $\{ q_n \}$ be its sequence of
  best returns. Then
  \begin{displaymath}
    \frac{1}{2 q_{n+1}} < \frac{1}{q_n + q_{n+1}} < \norm{ q_n
      \alpha} \leq \frac{1}{q_{n+1}}. 
  \end{displaymath}
\end{lem}
The second lemma shows that very good returns only occur for best
returns $q_n$ and their multiples.
\begin{lem}\label{lem:GoodReturns}
  Let $\alpha \in \R^+ \backslash \Q$ and let $\{ q_n \}$ be its sequence of
  best returns.If $q \in \Z \backslash \{0\}$ satisfies
  \begin{displaymath}
    \norm{q \alpha} < \frac{1}{2 |q|},
  \end{displaymath}
  then $q = l q_n$ for some best return $q_n$ and some
  \begin{displaymath}
    |l| < \sqrt{\frac{q_{n+1}}{q_n}}.
  \end{displaymath}
\end{lem}

\subsection{Cohomological Equations}

The behavior of the special flow $T^t_{(\alpha , \varphi)}$ is
determined by the cohomology class of the function $\varphi$,
see \cite{MR87f:28019}, \cite{Unpublished}, and \cite{MR1858535}.

We call two functions, $\varphi_1$ and $\varphi_2$, (additively)
cohomologous (over $R_\alpha$) if there is a measurable solution
$\psi$ to the equation 
\begin{displaymath}
  \psi (x + \alpha) - \psi(x) = \varphi_1(x) - \varphi_2(x).
\end{displaymath}
We call this equation the additive cohomological equation. If a
function is cohomologous to $0$ then we call it an (additive)
coboundary. Using this definition we can say $\varphi_1$ and
$\varphi_2$ are cohomologous if their difference $\varphi_1 -
\varphi_2$ is a coboundary. Coboundaries have mean $0$, hence no
positive function can be a coboundary. The appropriate notion of
triviality for positive functions is that of being cohomologous to a
constant. This corresponds to the associated special flow being
conjugate to a suspension flow.  Throughout our arguments we will use
that fact that we can subtract coboundaries from our function
$\varphi$ without altering the behavior of the special flow.
\begin{lem}
  Let $\varphi_1 : \T \rightarrow \R^+$ , $\varphi_2 : \T \rightarrow
  \R^+$, and $\alpha \in \R \backslash \Q$. If $\varphi_1$ and
  $\varphi_2$ are additively cohomologous, i.e. there exists a
  measurable ($L^2$) solution $\psi$ to the additive cohomological
  equation
  \begin{displaymath}
    \psi (x + \alpha) - \psi(x) = \varphi_1(x) - \varphi_2(x),
  \end{displaymath}
  then the special flow $T^t_{(\alpha,\varphi_1)}$ is measurably
  ($L^2$) conjugate to the special flow $T^t_{(\alpha, \varphi_2)}$.
  In particular, if $\varphi_2$ is a constant then the special flow
  $T^t_{(\alpha,\varphi)}$ is measurably ($L^2$) conjugate to a
  suspension flow.
\end{lem}

A cohomological equation again appears -- this time a multiplicative
cohomological equation -- when we study the existence of eigenvalues for the special flow.

\begin{lem}\label{lem:WeakMixingCriterion3}
  Let $\varphi : \T \rightarrow \R^+$, $\alpha \in \R \backslash \Q$,
  and $\lambda \in \R \backslash \{ 0 \}$. If there exists an
  increasing sequence $\{ m_n \}$ such that $\norm{ m_n \alpha }
  \rightarrow 0$ and
  \begin{displaymath}
    \int \norm{\lambda S_{m_n} \varphi} dx \not \rightarrow 0,
  \end{displaymath}
  then $\lambda$ is not an eigenvalue of the special flow
  $T^t_{(\alpha,\varphi)}$.
\end{lem}
\begin{proof}
  The eigenvalues of the special flow are determined by a
  multiplicative cohomological equation. In particular, $\lambda$ is
  an eigenvalue of the special flow if and only if there is a
  measurable solution $\Psi$ of the equation
  \begin{displaymath}
    e^{2 \pi i \lambda \varphi (x)} = \frac{\Psi(x + \alpha)}{\Psi(x)}.
  \end{displaymath}
  Iterating this we get for any $m$ the equation
  \begin{displaymath}
     e^{2 \pi i \lambda S_m \varphi (x)} = \frac{\Psi(x + m
     \alpha)}{\Psi(x)}
  \end{displaymath}
  and thus
  \begin{displaymath}
    e^{2 \pi i \lambda S_{m_n} \varphi (x)} -1 = \frac{\Psi(x + m_n
     \alpha)}{\Psi(x)}-1.
  \end{displaymath}
  By the property of $m_n$ that $\norm{ m_n \alpha} \rightarrow 0 $ we
  have that the right-hand side converges to $0$ in $L^1$. Thus the
  left-hand side also converges to $0$. By Lemma
  \ref{lem:EquivalentDistances}, if $\lambda$ is an eigenvalue, then
  \begin{displaymath}
    \int \norm{\lambda S_{m_n} \varphi} dx \rightarrow 0.
  \end{displaymath}
  Thus the given condition implies that $ \lambda$ is not an
  eigenvalue of the special flow. 
\end{proof}

The absence of eigenvalues other than the simple eigenvalue $0$, which
corresponds to the constant functions, implies weak mixing for the flow.
The eigenvalues for the flow form an additive subgroup of $\R$. Thus,
to prove the flow has a continuous spectrum, and is, hence, weak
mixing, it suffices to prove that no sufficiently large $\lambda$ is
an eigenvalue.

\subsection{Analytical Estimates}

We will analyze the cohomological equations via Fourier techniques. We
naturally arrive at considering expressions of the form $|e^{2 \pi i m
  \alpha}-1|$. These quantities are related to the quantities $\norm{m
  \alpha}$ appearing in Section \ref{sec:Arithmetic}. 
\begin{lem}
  \label{lem:EquivalentDistances}
  The two functions $\norm{x}$ and $|e^{2 \pi i x}-1|$ are related by
  \begin{displaymath}
    4 \norm{x} \leq |e^{2 \pi i x}-1| \leq 2 \pi \norm{x} .
  \end{displaymath}
\end{lem}

In order to use the criterion for the absence of an eigenvalue, Lemma
\ref{lem:WeakMixingCriterion3}, it is necessary to control the Birkhoff
sums of the function $\varphi$. When we consider these sums the
following lemma will be crucial. It is an immediate consequence of
Lemma \ref{lem:EquivalentDistances} that
\begin{equation}
  \label{eq:1}
  \frac{1}{2} \frac{ \norm{m k \alpha}}{\norm{ k \alpha}} < \Bigl|
  \frac{e^{2 \pi i m k \alpha} -1 } {e^{2 \pi i k \alpha}-1} \Bigr|
  < 2 \frac{\norm{m k \alpha}}{\norm{ k \alpha}}.
\end{equation}


\section{The Structure of the Proof}

We begin by looking for a conjugacy arising from the additive
cohomological equation. Supposing that the additive cohomological
equation has a solution given by a trigonometric sum and formally
solving for the necessary coefficients yields the formal series
\begin{displaymath}
  \psi (x) = \sum_{m \in \Z} \frac{c_m}{e^{2 \pi i m \alpha}-1}  e^{
  2\pi i m x}. 
\end{displaymath}
Using Lemma \ref{lem:EquivalentDistances} we estimate the coefficients
by
\begin{displaymath}
  \Bigl| \frac{c_m}{e^{2 \pi i m \alpha}-1} \Bigr| \leq \frac{|c_m|}{4
  \norm{m \alpha}}.
\end{displaymath}

If these coefficients are square summable, then the special flow
$T^t_{( \alpha ,\varphi)}$ is $L^2$ conjugate to a suspension flow.
Otherwise we must prove weak mixing. There are two different cases
depending on exactly how the sequence $\{ |c_m| / \norm{m \alpha} \}$
behaves.

\begin{prop}\label{prop:SingleFrequency}
  Let $\varphi$ be a $C^3$ function satisfying the hypotheses
  {\rm[H2]}, and {\rm[H3]}. If
  \begin{displaymath}
    \limsup_{m \rightarrow \infty} \frac{|c_m|}{\norm{m \alpha}} = C > 0,
  \end{displaymath}
  then the special flow $T^t_{(\alpha, \varphi)}$ is weak mixing.
\end{prop}

We call this the single frequency weak mixing case. In this case it
suffices to take a sequence $m_n$ in Lemma
\ref{lem:WeakMixingCriterion3} that consists of multiples of single
best returns as in \cite{Unpublished}. When this is not possible we use

\begin{prop}\label{prop:MultipleFrequency}
  Let $\varphi$ be a $C^3$ function satisfying the hypothesis
  {\rm[H1]}. If
  \begin{displaymath}
   \lim_{m \rightarrow \infty} \frac{|c_{m}|}{\norm{m \alpha}}=
  0 \qquad \text{and} \qquad \sum_{m \in \Z} \Bigl( \frac{|c_{m}|}{\norm{m
  \alpha}}\Bigr)^2 = \infty 
  \end{displaymath}
  then the special flow $T^t_{(\alpha, \varphi)}$ is weak mixing.
\end{prop}

We call this the multiple frequency weak mixing case. Our hypothesis
{\rm [H1]} ensures that the function $\varphi$ is cohomologous to a
function in which the only frequencies which appear are best returns.
In this case, it is not sufficient to take the $m_n$ to be a multiple
of a single frequency.  We will have to take $m_n$ to be a sum of many
frequencies. In this case, no frequency dominates, and, in fact, the
values of $S_{m_n} \varphi$ becomes normally distributed in the limit.

\section{Preliminary Reduction}

The cohomology classes of those functions $\varphi$ that appear in
Theorem \ref{thm:Dichotomy} admit nice representatives. We emphasize
that we get different ``well-adapted'' representatives for each
$\alpha \in \R \backslash \Q$.

\begin{lem}\label{lem:FirstReduction}
  Let $\varphi$ be a $C^3$ function. Let $\alpha \in \R \backslash \Q$
  have the sequence of best returns $\{ q_n \}$. Then $\varphi$ is
  cohomologous to the function $\varphi_1$ defined by
  \begin{displaymath}
    \varphi_1 (x) = \sum_{|m| \in M} c_m e^{2 \pi i m x},
  \end{displaymath}
  where $m \in M$ is either 0 or of the form $m = l q_n$, where $q_n$
  is a best return satisfying
  \begin{displaymath}
    q_{n+1} > q_n^2
  \end{displaymath}
  and $l$ is such that $l q_{n} < \frac{1}{2} q_{n+1}$.
\end{lem}

\begin{proof}
  Define the class $M$ by 
  \begin{equation}\label{eq:MDefinition}
    M := \{ m \geq 0: 2 m^2 \norm{m \alpha} \leq 1 \}
  \end{equation}
  and the function $\xi$ by
  \begin{displaymath}
    \xi (x) := \varphi(x) - \varphi_1(x) = \sum_{|m| \not \in M}
    c_m e^{2 \pi i m x}. 
  \end{displaymath}
  We need to show that $\xi$ is an additive coboundary. If $\xi (x) =
  \psi(x + \alpha ) - \psi(x)$ then $\psi(x)$ must be given by the
  formal series
  \begin{displaymath}
    \psi (x) := \sum_{|m| \not \in M} \frac{c_m}{e^{2 \pi i m
        \alpha}-1} e^{2 \pi i m x}.
  \end{displaymath}
  The coefficients of $\psi$ are estimated, using
  Lemma \ref{lem:EquivalentDistances} and \eqref{eq:MDefinition}, by
  \begin{displaymath}
    \Bigl| \frac{c_m}{e^{2 \pi i m \alpha}-1} \Bigr| \leq
    \frac{|c_m|}{4 \norm{ m \alpha}} \leq \frac{1}{2} m^2 |c_m| 
  \end{displaymath}
  Since $\varphi$ is $C^3$ these coefficients are square summable.
  Thus, the formal series is actually the Fourier series of an $L^2$
  function, and hence, $\xi$ is a coboundary.
  
  That $M$ contains only 0, best returns, and multiples of best
  returns follows from Lemma \ref{lem:GoodReturns}. The estimate on
  $q_{n+1}$ for $q_n \in M$ follows from Lemma
  \ref{lem:ContinuedFractionsEstimates}, and the definition of $M$
  \eqref{eq:MDefinition}, since
   \begin{displaymath}
     \frac{1}{2 q_{n+1}} < \norm{q_n \alpha} \leq \frac{1}{2 q_n^2}. 
   \end{displaymath}
\end{proof}

\section{Single Frequency Weak Mixing Case}

\subsection{Remarks}
Under hypothesis {\rm [H1]} the proof we give is strictly only
requisite for the case $C=\infty$. If $0 < C < \infty$, then
hypotheses {\rm [H2]} and {\rm [H3]} are not necessary since in this
case $\varphi$ is cohomologous to a function in which only best
returns appear, see Lemma \ref{lem:SecondReduction}. 

The argument for weak mixing given here is a classical one and is similar to that given by A.
Katok and E. A. Robinson in their 1983 unpublished notes
\cite{Unpublished}.


\subsection{Proof of Weak Mixing}
Under hypothesis {\rm [H2]} it is clear that $C$ must be achieved
along a subsequence $\{q_{s(n)} \}$ of best returns contained in $M$.
Let $\{ q_{s(n) } \}$ satisfy
\begin{displaymath}
  \lim_{n \rightarrow \infty} \frac{|c_{q_{s(n)}}|}{\norm{ q_{s(n)}
  \alpha}} = C.
\end{displaymath}

Fix $\lambda \in \R \backslash \{0\}$. We need to show that there
exists a sequence $m_n$ that serves in Lemma
\ref{lem:WeakMixingCriterion3} to show that $\lambda$ is not an
eigenvalue. For our sequence we take
\begin{equation}
  \label{eq:MNDefinition}
  m_n = b_n \qn := \Bigl\lceil \frac{\qnn}{4 \qn} \Bigr\rceil \qn.
\end{equation}
This sequence is chosen to isolate and inflate the terms corresponding
to multiples of the best return $\qn$. For this reason, it is natural
to consider the expression
\begin{equation}
  \label{eq:PhiNDefinition}
  \phi_n (x) = \sum_{|l| q_{s(n)} \in M} c_{l q_{s(n)}} e^{2 \pi i a
    q_{s(n)} x} .
\end{equation}

It is technically easier deal with a function that is nearly
constant. Since $\varphi$ is $C^3$ and any trigonometric polynomial
with $0$ average is an additive coboundary we may discard finitely
many terms from $M$ and suppose that
\begin{equation}
\label{eq:MPrimeDefinition}
\sum_{|m| \in M} \bigl| m c_m \bigr| < \frac{1}{16 |\lambda|}.
\end{equation}
We shall denote by $\varphi_\lambda$ the representative of the
cohomology class of $\varphi$ thus obtained.

We now show that the Birkhoff sums $S_{m_n} \varphi_\lambda$ are
uniformly close to the Birkhoff sums $S_{m_n} \phi_n$ for all $n$. The
Birkhoff sums $S_{m_n} \phi_n$ are much simpler and we will be able to
estimate them directly.

\begin{lem}\label{lem:BirkhoffSumEstimate1}
Suppose $\varphi_\lambda$ satisfies \eqref{eq:MPrimeDefinition} and
$\phi_n$ is given by \eqref{eq:PhiNDefinition}.  For all $n$,
\begin{equation}
  \label{eq:BirkhoffSumEstimate1}
  \bigl| \lambda S_{m_n} \varphi_\lambda (x) - \lambda S_{m_n} \phi_n
  (x) \bigr| < \frac{1}{8}.
\end{equation}
\end{lem}

\begin{proof}
We will show that
\begin{displaymath}
  \bigl| \lambda S_{m_n} \varphi_\lambda (x) - \lambda S_{m_n} \phi_n
  (x) \bigr|   \leq 2 | \lambda | \sum_{|m| \in M} | m c_m |
\end{displaymath}
from which we get the required estimate using
\eqref{eq:MPrimeDefinition}. We can directly compute
\begin{displaymath}
  \bigl| \lambda S_{m_n} \varphi_\lambda (x) - \lambda S_{m_n} \phi_n
  (x) \bigr| \leq |\lambda| \sum_{|l| q_k \in M: k \neq
  s(n)} \Bigl| \frac{e^{2 \pi i m_n l q_k \alpha}-1}{e^{2 \pi i l q_k
  \alpha}-1} \Bigr| |c_{l q_k}|.
\end{displaymath}
Using the estimate \eqref{eq:1} yields
\begin{equation}
  \label{eq:BasicEstimate1}
  \Bigl| \frac{e^{2 \pi i m_n l q_k \alpha}-1}{e^{2 \pi i l q_k
  \alpha}-1} \Bigr| \leq 2 \frac{\norm{m_n q_k
  \alpha}}{\norm{q_k \alpha}} \leq 2 m_n.
\end{equation}
For $k > s(n)$ we have $2 m_n < |l q_k|$. For $k < s(n)$ we need to use
the fact that $m_n$, as a multiple of $\qn$, produces a better return
than does $q_k$. From the definition of $m_n$ \eqref{eq:MNDefinition}
we get
\begin{displaymath}
  \frac{\norm{ m_n q_k \alpha}}{\norm{ q_k \alpha}} \leq \frac{b_n q_k
    \norm{\qn \alpha} }{ \norm{q_k \alpha}}.
\end{displaymath}
Using Lemma \ref{lem:ContinuedFractionsEstimates} and estimating $b_n
< \frac{\qnn}{2\qn}$ yields
\begin{displaymath}
  \frac{\norm{ m_n q_k \alpha}}{\norm{ q_k \alpha}} < \frac{ q_k
    q_{k+1}}{\qn} < q_k.
\end{displaymath}
For all $|l q_k| \in M$ with $k \neq s(n)$ we have
\begin{displaymath}
  \Bigl| \frac{e^{2 \pi i m_n l q_k \alpha}-1}{e^{2 \pi i l q_k
  \alpha}-1} \Bigr| \leq 2  \frac{\norm{m_n q_k
  \alpha}}{\norm{q_k \alpha}}  \leq |l| q_k,
\end{displaymath}
which proves the result.
\end{proof}

We estimate the Birkhoff sums $S_{m_n} \phi_n$ geometrically. In
essence, we show that an appropriately renormalized version of the sum
has a derivative that is of the same magnitude as the length of the
range. Since the length of the range does not go to zero, this is
sufficient to conclude that the integral in Lemma
\ref{lem:WeakMixingCriterion3} does not go to zero.

\begin{lem}\label{lem:BirkhoffSumEstimate2}
  Let $\varphi_n$, given by \eqref{eq:PhiNDefinition} satisfy
  hypotheses {\rm [H2]} and {\rm [H3]}. For $\lambda > 0$ sufficiently
  large there exists $\epsilon > 0$ such that
\begin{displaymath}
  \mu \bigl\{ x : \norm{ \lambda S_{m_n} \phi_n (x) } \geq \frac{1}{4}
  \bigr\} > \epsilon 
\end{displaymath}
for all $n$ sufficiently large.
\end{lem}

\begin{proof}
  Let
  \begin{align*}
    R_n &:= \sup \,  S_{m_n} \phi_n (x) - \inf \,  S_{m_n}
    \phi_n (x)\\
     D_n  &:= \sup |  S_{m_n} \phi_n' (x)|
  \end{align*}
  We can estimate $R_n$ from below using [H2] as
  \begin{align*}
    \label{eq:RangeEstimate}
    R_n &> 4 \Bigl( \Bigl| \frac{e^{2 \pi i m_n \qn \alpha}-1}{e^{2
        \pi i \qn \alpha}-1} c_{\qn}\Bigr| - \sum_{l=2}^{\infty}\Bigl|
    \frac{e^{2 \pi i m_n l \qn \alpha}-1}{e^{2 \pi i l \qn \alpha}-1}
    c_{l \qn}\Bigr| \Bigr)
    \\
    &> 4 \frac{\norm{m_n \qn \alpha}}{\norm{\qn \alpha}} \Bigl(
    \frac{1}{2} |c_{\qn}| - 2 \sum_{l=2}^{\infty} |c_{l \qn}| \Bigr)
    \\
    &> 2 \frac{\norm{m_n \qn \alpha}}{\norm{\qn \alpha}}\, |c_{\qn}|
        (1 - 4 K_1).
  \end{align*}
  Similarly, we can estimate $D_n$ from above using [H3] as
  \begin{align*}
    D_n & < 2\, \sum_{l=1}^{\infty}\Bigl| \frac{e^{2 \pi i m_n l \qn
        \alpha}-1}{e^{2 \pi i l \qn \alpha}-1} \,l \qn \, c_{a
      \qn}\Bigr| \\
    & < 4 \, \qn \, \frac{\norm{m_n \qn \alpha}}{\norm{\qn \alpha}} \,
    |c_{\qn}| (1+ K_2) .
  \end{align*}
  Now consider the number $I$ of intervals of the form $\bigl[ p +
  \frac{1}{4}, p + \frac{3}{4} \bigr]$ contained in the range of
  $\lambda S_{m_n} \phi_n$. This can be estimated from below by $I >
  \lfloor |\lambda| R_n \rfloor -1$. For $\lambda$ and $n$
  sufficiently large, $|\lambda| R_n > 4$ and hence $I > |\lambda| R_n
  /2 $. By continuity, $\lambda S_{m_n} \phi_n$ must cross each
  interval at least once and hence, by periodicity, it must cross each
  interval at least $\qn$ times.  Therefore, we have for each interval
  contained in the range of $\lambda S_{m_n} \phi_n$
  \begin{displaymath}
    \mu \Bigl\{x :  \lambda S_{m_n} \phi_n (x) \in \bigl[ p +
    \frac{1}{4}, p + \frac{3}{4} \bigr] \Bigr\} > \frac{\qn}{2
    |\lambda| D_n}.
  \end{displaymath}
  Multiplying this by our estimate for the number of intervals we get
  \begin{displaymath}
    \mu \bigl\{ x : \norm{ \lambda S_{m_n} \phi_n (x) } \geq \frac{1}{4}
    \bigr\}  >  \frac{\qn R_n}{4 D_n} >  \frac{ 1- 4 K_1}{8 (1+ K_2)} >0.
  \end{displaymath}
\end{proof}

Thus, combining Lemma \ref{lem:BirkhoffSumEstimate1} and Lemma
\ref{lem:BirkhoffSumEstimate2} we have, for $\lambda$ sufficiently
large, that
\begin{displaymath}
  \mu \bigl\{ x : \norm{ \lambda S_{m_n} \varphi_\lambda (x) } \geq
  \frac{1}{8} \bigr\} >  \epsilon
\end{displaymath}
for all $n$. Using our criterion for the absence of an eigenvalue,
Lemma \ref{lem:WeakMixingCriterion3}, this shows that $\lambda$ is not an
eigenvalue of the special flow. Since $\lambda$ was any sufficiently
large number this shows, by the remark following Lemma
\ref{lem:WeakMixingCriterion3}, that the special flow is weak mixing.

\section{Multiple Frequency Weak Mixing Case}

We simplify our problem by extracting an even simpler representative
of the cohomology class of $\varphi$. At this point the multiples of
the best returns are used. It is at this point that we use hypothesis
[H1].

\begin{lem}\label{lem:SecondReduction}
  Let $\varphi$ be a $C^3$ function satisfying hypothesis {\rm[H1]}. If
  \begin{displaymath}
    \sup \frac{|c_m|}{\norm{ m \alpha}}  = K_3 < \infty,
  \end{displaymath}
  then the function $\varphi$ is cohomologous to the function
  $\varphi_2$ given by
  \begin{displaymath}
    \varphi_2 (x) = \sum_{|m| \in M'} c_m e^{2 \pi i m x},
  \end{displaymath}
  where $m \in M'$ is either $0$ or a best return $q_n$ satisfying
  $q_{n+1} >q_n^2$.
\end{lem}

\begin{proof}
  Applying Lemma \ref{lem:FirstReduction} we see it suffices to prove
  that we may exclude the multiples of best returns. Let $\xi$ be the
  trigonometric series generated by the multiples of best returns,
  \begin{displaymath}
    \xi (x) = \sum_{|l| q_n \in M : |l| \geq 2} c_{l q_n}
    e^{ 2 \pi i l q_n x}.
  \end{displaymath}
  If $\xi(x) = \psi(x + \alpha) -  \psi(x)$, then $\psi$ must be
  given by the formal series
  \begin{displaymath}
    \psi(x) = \sum_{|l| q_n \in M: |l| \geq 2}
    \frac{c_{l q_n}}{e^{2 \pi i l q_n \alpha}-1}
    e^{ 2 \pi i l q_n x}.
  \end{displaymath}
  Using [H1] we get
  \begin{displaymath}
    \sum_{a = 2}^\infty \Bigl| \frac{c_{l q_n}}{e^{2 \pi i l q_n
        \alpha}-1} \Bigr|^2  \leq \frac{1}{ 16 \norm{q_n \alpha}^2}
    \sum_{l=2}^\infty \frac{|c_{l q_n}|^2}{ a^2 } < \frac{C_{q_n}}{16} K_3^2.
  \end{displaymath}
  Thus, since by [H1] the $C_m$ are summable, we have have $\psi$ is
  actually an $L^2$ function and hence $\xi$ is an additive coboundary.
\end{proof}

No individual frequency contributes enough to prove weak mixing. In
order to apply our criterion we need to take a group of frequencies
together. Our weak mixing sequence is of the form
\begin{displaymath}
  m_n := \sum_{k = l_n}^{u_n} b_k q_{s(k)},
\end{displaymath}
where 
\begin{displaymath}
  b_k := \Bigl\lceil \frac{q_{s(k)+1}}{4 q_{s(k)}} \Bigr\rceil
\end{displaymath}
and $l_n$ is an increasing sequence. This satisfies the requirement
that $\norm{m_n \alpha} \rightarrow 0$ regardless of the exact choices
of $l_n$ and $u_n$. 

We now prove a lemma analogous to Lemma \ref{lem:BirkhoffSumEstimate1}
for our more complicated situation. Our sequence $m_n$ is chosen to
isolate and inflate those coefficients corresponding to the
frequencies $\{ q_{s(k)} \}_{k =l_n}^{u_n}$. For this reason it is
natural to define
\begin{equation}
  \label{eq:ChiNDefinition2}
  \phi_n (x) :=c_{-\qn} e^{-2 \pi i \qn x} + c_0 + c_{\qn} e^{2 \pi i \qn
  x}.
\end{equation}
These functions asymptotically capture all the behavior in $S_{m_n}
\varphi (x)$. Unfortunately, their behavior is not as easy to control
as in the single frequency case.
\begin{lem}\label{lem:BirkhoffSumEstimate3}
  Let
  \begin{equation}
    \label{eq:MainDifference}
    \Delta_n = \bigl\| S_{m_n} \varphi (x) - 
    \sum_{k=l_n}^{u_n} S_{b_k q_{s(k)}} \phi_k \bigl(x +
    \sum_{j=l_n}^{k-1} b_j q_{s(j)} \alpha \bigr) \bigr\|_{\infty}.
  \end{equation}
  Then $\Delta_n \rightarrow 0$ as $n \rightarrow \infty$.
\end{lem}

\begin{proof}
  Define
  \begin{displaymath}
    \delta_k := \bigl\| S_{b_k q_{s(k)}} \varphi (x) - S_{b_k q_{s(k)}} \phi_k (x) \bigr\|_{\infty}
  \end{displaymath}
  and observe that
  \begin{displaymath}
    \Delta_n \leq \sum_{k = l_n}^{u_n} \delta_k .
  \end{displaymath}
  If we show that $\delta_k$ is a summable sequence then $\lim_{n
    \rightarrow \infty} \Delta_n = 0$ follows from $l_n \rightarrow
  \infty$ and is independent of the choice of $u_n$. Using the
  triangle inequality and Lemma \ref{lem:EquivalentDistances} we get
  \begin{displaymath}
    \delta_k \leq 4 \sum_{j \in \N \backslash \{k\} }
    \frac{\norm{q_{s(j)} b_k q_{s(k)} \alpha}} {\norm{q_{s(j)}
        \alpha}} \bigl| c_{q_{s(j)}}\bigr| := 4 \sum_{j \in \N
    \backslash \{k\} } Q_j 
  \end{displaymath}
  We break the sum into two pieces, which we will estimate separately,
  \begin{displaymath}
    \delta_k \leq 4 \sum_{j=1}^{k-1} Q_j + 4 \sum_{j=k+1}^{\infty}
    Q_j 
  \end{displaymath}
  For $j > k$, using the bound for $|c_m| / \norm{m \alpha}$, we
  produce
  \begin{displaymath}
    Q_j = b_k q_{s(k)} \norm{
    q_{s(j)} \alpha} \frac{|c_{q_{s(j)}}|}{\norm{q_{s(j)} \alpha}}
    \leq \frac{K q_{s(k)+1}}{2 q_{s(j)+1}}.
  \end{displaymath}
  Using the condition $q_{s(k)+1} > q_{s(k)}^2$ yields
  \begin{displaymath}
      \sum_{j=k+1}^{\infty} \frac{K q_{s(k)+1}}{2 q_{s(j)+1}} \leq
      \sum_{j= k+1}^{\infty} \frac{K q_{s(k+1)}}{2 q_{s(j)}^2} \leq
      \frac{K}{q_{s(k+1)}^2}. 
  \end{displaymath}
  This is summable in $k$. For $j < k$ we produce an estimate
  analogous to the one above using the fact that $q_k$ produces a
  better return than does $q_j$, 
  \begin{displaymath}
    Q_j = b_k q_{s(j)} \norm{q_{s(k)} \alpha} \frac{|c_{q_{s(j)}}|}
    {\norm{q_{s(j)} \alpha}} \leq \frac{K_3 q_{s(j)}}{2 q_{s(k)}} .
  \end{displaymath}
  Again, using the condition $q_{s(k)+1} > q_{s(k)}^2$,
  we obtain
  \begin{displaymath}
    \sum_{j=1}^{k-1} \frac{K_3 q_{s(j)}}{2 q_{s(k)}} \leq \frac{K_3 k
    q_{s(k-1)}}{2 q_{s(k)}} \leq \frac{K_3 k}{2 q_{s(k-1)}}    
  \end{displaymath}
  which is summable in $k$. 
\end{proof}

Thus, $S_{m_n} \varphi (x)$ is asymptotic to $\sum_{k = l_n}^{u_n}
S_{b_k q_{s(k)}} \phi_k (x)$. Ignoring the constant term, this sum is
of the form
\begin{displaymath}
  \sum_{k = l_n}^{u_n} S_{b_k q_{s(k)}} \phi_k (x) =
  \sum_{k=l_n}^{u_n} d_k \cos ( q_{k} x + r_k ).
\end{displaymath}
We choose a sequence $(l_n,u_n)$ with $l_n \rightarrow \infty$ and
such that
\begin{equation}\label{eq:Variance}
  \lim_{n \rightarrow \infty} \sum_{k=l_n}^{u_n} d_k^2 = 1.
\end{equation}
For lacunary series $q_n$, the random variables $\cos(q_n x + r_n)$ are
only weakly dependent. This observation will allow us to compute the
asymptotic distribution of the sums $S_{m_n} \phi_n (x)$.

\subsection{The Distribution of the Birkhoff Sums}

The proof of the usual central limit theorem for lacunary
trigonometric sum was carried out by Salem and Zygmund in 1947
\cite{MR9:181d}. Unfortunately, the convergence of normalized sums is
not exactly what is needed. What is needed is a version of the
``series'' central limit theorem \cite{MR85a:60007} for lacunary series.
Fortunately, the proof of Salem and Zygmund carries through with no
changes to prove this theorem.

\begin{thm}
  Let
  \begin{displaymath}
    X_n (x) = \sum_{k = 1}^{u_n} c_{k,n} \cos(q_{k,n} x + r_{k,n})
  \end{displaymath}
  for some sequence $q_{k+1,n} \geq \lambda q_{k,n}$ with $\lambda >
  1$ and some coefficients satisfying
  \begin{displaymath}
    \var(X_n (x)) = \sum_{k=1}^{u_n} c_{k,n}^2 = 1
  \end{displaymath}
  and 
  \begin{displaymath}
    c_{k,n} \rightarrow 0 \text{ uniformly as $n \rightarrow \infty$}.
  \end{displaymath}
  Then 
  \begin{displaymath}
    X_n \xrightarrow{dist} N(0,1),
  \end{displaymath}
  where $N(0,1)$ is the normal distribution with mean $0$ and variance $1$.
\end{thm}

\begin{proof}
We use the method of characteristic functions developed by Lyapunov to
prove the central limit theorem. For simplicity we only prove the case
where $q_{k+1,n} > 2 q_{k,n}$ since this is sufficient for us.

Let $F_n$ denote the distribution functions for the random variables
$X_n$. Let $\varphi_n (t)$ be the characteristic function of the
distribution $F_n$,
\begin{displaymath}
  \varphi_n (t) = \int_{-\infty}^{+\infty} e^{ i t y} d F_n (y) .
\end{displaymath}
Our goal is to show that these characteristic functions converge to
that of the normal distribution,
\begin{displaymath}
  \lim_{n \rightarrow \infty} \varphi_n (t) = e^{ - \frac{t^2}{2}}.
\end{displaymath}

Passing from the integral with respect to the distribution to the
integral with the random variable $X_n$ we get
\begin{displaymath}
  \varphi_n (t) = \frac{1}{2 \pi} \int_{0}^{2 \pi} \exp \bigl( i t
  \sum_{k = 1}^{u_n}  \cos(q_{k,n} x + r_{k,n}) \bigr) dx .
\end{displaymath}
We use the relation
\begin{displaymath}
  e^z = (1 + z ) e^{\frac{1}{2} z^2 +o(|z^2|)} 
\end{displaymath}
and the fact $\sum_{k=1}^{u_n} c_{k,n}^2 = 1$ to obtain
\begin{displaymath}
  \begin{split}
    \varphi_n (t) 
    = \frac{1}{2 \pi} \int_{0}^{2 \pi} e^{o(1)}
  \prod_{k=1}^{u_n}\bigl( 1 +& i t c_{k,n} \cos (q_{k,n} x + r_{k,n})
  \bigr)\\
  &\exp \bigl( - \frac{1}{2} t^2 c_{k,n}^2 \cos (q_{k,n} x +
  r_{k,n}) \bigr) dx.    
  \end{split}
\end{displaymath}

First we show that the first term is bounded. 
\begin{displaymath}
  \Bigl| \prod_{k = 1}^{u_n} \bigl( 1+ i t c_{k,n} \cos (q_{k,n} x +
  r_{k,n}) \bigr)  \Bigr| \leq \prod_{k= 1}^{u_n} \bigl( 1 + t^2 c_{k,n}^2
  \bigr)^{\frac{1}{2}} \leq e^{\lambda^2} 
\end{displaymath}
The exponent of the second term can be rewritten using the double
angle formula as 
\begin{displaymath}
  \sum_{k=1}^{u_n} c_{k,n}^2 \cos^2 (q_{k,n} x + r_{k,n}) = 1 +
 \sum_{k=1}^{u_n}  \frac{ c_{k,n}^2}{2} \cos (2 q_{k,n} x + 2
 r_{k,n} ) = 1+ \xi_n (x). 
\end{displaymath}
The Lebesgue measure of the set of points where $|\xi_n (x) | \geq
\delta >0$ can be estimated by 
\begin{displaymath}
  \frac{1}{\delta^2} \int_{0}^{2 \pi} \xi_n^2 (x) = \frac{1}{4
  \delta^2} \pi \sum_{k=1}^{u_n} c_{k,n}^4 \stackrel{N \rightarrow
  \infty}{\rightarrow } 0    
\end{displaymath}
since $c_{k,n} \rightarrow 0$. Since $\xi_n$ is bounded we have that the
convergence in measure implies $L^1$ convergence from which we get
that the integral is asymptotic to 
\begin{displaymath}
  e^{-\frac{t^2}{2}} \int_{0}^{2 \pi} \prod_{k=1}^{u_n} \bigl( 1+ i 
  t c_{k,n} \cos (q_{k,n} x + r_{k,n}) \bigr) dx.
\end{displaymath}
Since the first term is the characteristic function of the Gaussian
random variable with mean $0$ and variance $1$ we simply need to show
that 
\begin{equation}\label{eq:CosineProduct}
  \int_{0}^{2 \pi} \prod_{k=1}^{u_n} \bigl( 1+ i t c_{k,n}
  \cos (q_{k,n} x + r_{k,n}) \bigr) dx = 1.
\end{equation}
Using the fact that
\begin{displaymath}
  \cos (m x ) \cos (n x) = \frac{1}{2} \bigl( \cos \bigl( (m+n) x\bigr) +
  \cos \bigl( (m-n) x \bigl) \bigr),
\end{displaymath}
we can rewrite the product in the integral in the form
\begin{displaymath}
  \prod_{k=1}^{u_n} \bigl( 1+ i t c_{k,n}
  \cos (q_{k,n} x) \bigr) = a_{0,n} + \sum_{l \in L_n} a_{l,n} \cos (l
  x+ s_l) ,
\end{displaymath}
where
\begin{displaymath}
  L_n := \bigl\{ l \geq 0 \, : \, l = \sum_{k = 1}^{u_n} b_{k} q_{k,n}
  \text{ with } b_k \in \{ -1 , 0 ,1\} \bigr\}. 
\end{displaymath}
Given that $q_{k+1,n} \geq 2 q_{k,n}$, we immediately observe
$\sum_{k=1} {l-1}q_{k,n} < q_{l,n}$ from which it immediately follows
that the representation of $0$ in the form $\sum b_k q_{k,n}$ is
impossible.  Thus we have $a_{0,n} =1$. Integrating, we immediately get
the proof of the theorem.
\end{proof}

\subsection{Proof of Weak Mixing}

Let $z_n = m_n \lambda c_0 \mod 1$. By passing to a subsequence we may assume
that $ z_n \rightarrow z $. Using \eqref{eq:Variance}, yields
\begin{displaymath}
  \lambda S_{m_n} \phi_n \xrightarrow{d} N(z,1).
\end{displaymath}
Finally,
\begin{displaymath}
  \lim_{n \rightarrow \infty} \int \norm{\lambda S_{m_n} \varphi (x) }
  dx = \lim_{n \rightarrow \infty} \int \norm{ z_n + \lambda S_{m_n}
  \phi_n (x) } = \int \norm{x} dw(x) > 0,
\end{displaymath}
where $w$ is the measure corresponding to the $N(z,1)$ distribution.
This proves that $ \lambda$ is not an eigenvalue. Since $\lambda \in \R
\backslash \{ 0 \}$ was arbitrary this proves that the special flow
has continuous spectrum and is hence weak mixing.
 

\section{Final Comments}

Almost all of our arguments require very weak hypotheses. Every
statement with the exception of Lemma \ref{lem:SecondReduction} holds
for $C^3$ functions satisfying the appropriate regularity of decay
properties. We should be able to considerably enlarge the class of
functions for which our dichotomy holds by using appropriate uniform
estimates similar to Lemma \ref{lem:BirkhoffSumEstimate1} rather than
the stronger statements of Lemma \ref{lem:SecondReduction} and Lemma
\ref{lem:BirkhoffSumEstimate3}.

\bibliographystyle{plain}
\bibliography{Dichotomy.bib}

\end{document}